\documentclass[12pt]{article}

\usepackage{amsmath}
\usepackage{lipsum}
\usepackage{tikz-cd}

\usepackage{amsthm}
\usepackage{mathrsfs}	
\usepackage{amssymb}
\usepackage{amscd}
\usepackage{setspace}	
\usepackage{tikz}		
\usepackage{tikz-cd}
\usetikzlibrary{matrix}	
\usepackage{graphicx}	
\graphicspath{ {pics/} }	
\usepackage{wrapfig}	
\usepackage{enumerate}	
\usepackage{bbm}		

\usepackage[utf8]{inputenc}	
\usepackage[english]{babel}	
\usepackage{blindtext}		
\usepackage[affil-it]{authblk}	

\usepackage{faktor}			
\usepackage{xparse}			

\usepackage{hyperref}		
\usepackage[lite]{amsrefs}

\newenvironment{acknowledgements}{%
  \begin{abstract}
}{%
  \end{abstract}
}

\newtheorem{theorem}{Theorem}[section] 

\newtheorem{definition}[theorem]{Definition}

\newtheorem{lemma}[theorem]{Lemma}

\newtheorem{remark}[theorem]{Remark}

\newtheorem{proposition}[theorem]{Proposition}

\newtheorem{corollary}[theorem]{Corollary}

\newtheorem{examples}[theorem]{Examples}

\newtheorem{notation}[theorem]{Notation}

\newcommand{\acts}{\curvearrowright}
\DeclareMathOperator{\asdim}{\text{asdim}}
\newcommand{\diam}{\text{diam}}

\newcommand{\Z}{\mathbb{Z}}

\newcommand{\dad}{\text{DAD}}

\newtheorem{thmx}{Theorem} 

\newtheorem*{theorem*}{Example}

\DeclareFontFamily{U}{mathb}{\hyphenchar\font45}
\DeclareFontShape{U}{mathb}{m}{n}{ <-6> matha5 <6-7> matha6 <7-8>
mathb7 <8-9> mathb8 <9-10> mathb9 <10-12> mathb10 <12-> mathb12 }{}
\DeclareSymbolFont{mathb}{U}{mathb}{m}{n}

\DeclareMathAccent{\abxring}{0}{mathb}{"38}

\DeclareFontFamily{U}{mathb}{\hyphenchar\font45}
\DeclareFontShape{U}{mathb}{m}{n}{ <-6> matha5 <6-7> matha6 <7-8>
mathb7 <8-9> mathb8 <9-10> mathb9 <10-12> mathb10 <12-> mathb12 }{}
\DeclareSymbolFont{mathb}{U}{mathb}{m}{n}

\begin{document}
\title{A Hurewicz-type Theorem for the Dynamic Asymptotic Dimension with Applications to Coarse Geometry and Dynamics}
\author{Samantha Pilgrim}
\maketitle

\paragraph{Abstract:} We prove a Hurewicz-type theorem for the dynamic asymptotic dimension originally introduced by Guentner, Willett, and Yu.  Calculations of (or simply upper bounds on) this dimension are known to have implications related to cohomology of group actions and the $K$-theory of their transformation group $C^*$-algebras.  Moreover, these implications are relevant to the current classification program for $C^*$-algebras.  As a corollary of our main theorem, we show the dynamic asymptotic dimension of actions by groups on profinite completions along sequential filtrations by normal subgroups is often subadditive over extensions of groups, which shows that many such actions by elementary amenable groups are finite dimensional.  We combine this extension theorem with other novel results relating the dynamic asymptotic dimension of such actions to the asymptotic dimension of corresponding box spaces.  This allows us to give upper bounds on the asymptotic dimension of many box spaces (including those of infinitely-many groups with exponential growth).  For some of these examples, we can also find lower bounds by utilizing the theory of ends of groups.  

\section{Introduction}
When considering any invariant thought to behave like dimension, a natural question is whether it is subadditive over morphisms.  That is, given a morphism $f: X\to Y$ of two objects, is the dimension of $X$ at most the dimension of $Y$ plus the largest dimension of a `fibre' of $f$  (speaking loosely)?  An affirmative answer to such a question is sometimes called a Hurewicz-type theorem, and speaks to the naturality and dimension-theoretic nature of such an invariant.  Classical results of this kind include the rank-nullity theorem (for the category of vector spaces and linear maps) and the Hurewicz mapping theorem \cite[Theorem 1.12.4, page 136]{engelking} (for the category of topological spaces and continuous maps).  Special cases include whether a notion of dimension is subadditive over products or extensions of objects.  

More recently, Hurewicz-type theorems were developed for Gromov's asymptotic dimension in \cite{asymptotichurewicztheorem} and later in \cite{Brodskiy2006AHT}.  This paper uses similar methods as in \cite{Brodskiy2006AHT}, together with some original ideas, to prove a Hurewicz-type theorem for the dynamic asymptotic dimension ($\dad$) of Guentner, Willett, and Yu.  This dimension theory was introduced in \cite{dasdimGWY} as a dynamical analogue of the asymptotic dimension.  It is related to conditions used in work on the Farrell-Jones conjecture on manifold topology \cite{bartels} \cite{bartels_et_al}, as well as to the homology of group actions \cite{bonicke2021dynamic} and the nuclear dimension of transformation group $C^*$-algebras \cite[Theorem 8.6]{dasdimGWY}.  It has also been applied to the hyperfiniteness problem for orbit equivalence relations \cite{Conley2020BorelAD}.  

We state the appropriate Hurewicz-type theorem for the dynamic asymptotic dimension below.  We use the definition of $\dad$ for free actions (see \ref{dad definition}), in the statement, but do not assume the actions are free (more on this in \ref{freeness remark}).  

\begin{thmx}
Suppose $\Gamma\acts X$ and $\Lambda\acts Y$ are actions on compact metric spaces by homeomorphisms, and that $f: X\to Y$ and $\alpha: \Gamma\times X\to \Lambda$ are continuous functions such that $f(\gamma\cdot x) = \alpha(\gamma, x)\cdot f(x)$ for all $x\in X$ and $\gamma\in \Gamma$.  If $\dad(f):=\sup_{A\subset Y : \dad(A) = 0} \dad_{free}(f^{-1}(A))$, it holds that 

$$\dad_{free}(\Gamma\acts X)\leq \dad(f) + \dad_{free}(\Lambda\acts Y).$$  
\end{thmx}\label{Theorem A}

Apart from being a natural permanence property for the dynamic asymptotic dimension, this theorem is interesting because it provides a new tool for showing certain actions have finite $\dad$.  While it is thought that an action by an amenable group with finite asymptotic dimension should have finite $\dad$, relatively little progress has been made aside from the case where the acting group is virtually nilpotent \cite[Theorem 7.4]{szabo2017rokhlin}.  Finiteness of the dynamic asymptotic dimension is important to the current classification program for $C^*$-algebras because of how $\dad$ relates to nuclear dimension and computations of $K$-theory via the Baum-Connes and HK conjectures \cite{Guentner2016DynamicAD} \cite{bonicke2021dynamic}.  New examples of finite dimensional actions add to the class of algebras which are both classifiable by their $K$-theoretic and tracial data and which admit extensive technology for computing $K$-theory.  $\dad$ also plays a special role for those approaching classification from a dynamical perspective, as it applies to actions with finite stabilizers which are not essentially free, and therefore not covered by Kerr's notion of almost finitness (see for instance \cite{MR4066584} and \cite{https://doi.org/10.48550/arxiv.2107.05273}).  The introduction of \cite{https://doi.org/10.48550/arxiv.2201.03409} gives a precise statement of the classification theorem and additional references.  

We give two applications of this new Hurewicz-type theorem.  First, we can show that the $\dad$ is preserved when one restricts to a finite-index subgroup of the acting group, similar to how the asymptotic dimension of a finite index subgroup is the same as that of the larger group.  This further demonstrates that the $\dad$ is a coarse dimension theory for dynamics.  

\begin{thmx} Suppose $\Gamma\acts X$ is a free action on a compact metric space by homeomorphisms and $\Lambda<\Gamma$ is a finite index subgroup.  Then $\dad(\Gamma\acts X)= \dad(\Lambda\acts X)$.  
\end{thmx}\label{Theorem B}

\noindent We can also show that the $\dad$ of groups acting on their profinite completions (sometimes called odometer actions) behaves subadditively when one has an extension of groups, at least when the extension has certain properties:

\begin{thmx}
Suppose $1\to \Delta\to \Gamma\to \Lambda\to 1$ is an exact sequence of countable groups and that $\Gamma$ is residually finite.  Suppose $(N_i)$ is a countable collection of finite index normal subgroups of $\Gamma$ which separates $\Delta/(\cap_i (N_i\cap \Delta))$ from $\Gamma/(\cap_i N_i)$ (see \ref{separable subgroup definition}).  This collection induces collections of finite index normal subgroups and hence odometers of $\Delta$ and $\Lambda$.  Then if $\text{DAD}_{free}(\Delta\acts \widehat{\Delta}_{(N_i\cap \Delta)}) = m$ and $\text{DAD}_{free}(\Lambda\acts \widehat{\Lambda}_{(N_i/\Delta\cap N_i)}) = n$, $\text{DAD}_{free}(\Gamma\acts \widehat{\Gamma}_{(N_i)}) \leq m + n$.  
\end{thmx}\label{Theorem C}

This extension theorem can further be applied to coarse geometry, specifically the study of box spaces of groups.  These are spaces constructed out of finite quotients of a group.  Just as the asymptotic dimension of a group itself can have interesting group-theoretic applications (see for instance \cite{MR2431020}), there are likewise many correspondences between geometric properties of box spaces and group-theoretic properties.  For example, coarse equivalence of box spaces encodes virtual isomorphism of groups \cite{MR4091056}, and amenability of groups is equivalent to property A of their box spaces \cite[Proposition 11.39]{Roe2003LecturesOC}.  Similar results also exist for a-T-menablity \cite{https://doi.org/10.1112/blms/bdt045}.  

Box spaces provided some of the first examples of spaces which do not coarsely embed in Hilbert space \cite[Remark 11.27]{Roe2003LecturesOC} and it has been shown that a group with a finite dimensional box space must be amenable \cite[Remark 11.38 and Proposition 11.39]{Roe2003LecturesOC}.  Just as actions by finite dimensional, amenable groups are thought to be finite dimensional, it is thought that a partial converse to this latter fact should also hold; that is, an amenable group with finite asymptotic dimension ought to have finite dimensional box spaces.  Notwithstanding, there has been similarly little progress on this aside from the case of virtually nilpotent groups \cite{boxspacesDT}.  

As it turns out, questions about the asymptotic dimension of box spaces are essentially special cases of questions about the dynamic asymptotic dimension, echoing a relationship between DAD and the asymptotic dimension of warped cones described in \cite[Section 8]{warpedcones}.  More precisely, we prove the following in \ref{odometers and box spaces}.  

\begin{thmx} \label{Theorem D}
Let $\Gamma$ be a countable group and $(N_i)$ a sequence of finite index subgroups of $\Gamma$.  Suppose $\Gamma\acts \widehat{\Gamma}_{(N_i)}$ is the action by left multiplication of $\Gamma$ on its profinite completion along $(N_i)$.  Then $\dad_{free}(\Gamma\acts \widehat{\Gamma}_{(N_i)}) = \asdim\Box_{(N_i)}\Gamma$.  
\end{thmx}

This relationship allows us to prove a special case of an extension theorem for the asymptotic dimension of box spaces.  In fact, such an extension theorem was our original motivation for investigating a Hurewicz-type theorem for the DAD.  A proof of a similar result was attempted in \cite{finnsell2015asymptotic}, but the proof therein contains a gap which we discuss later in this section.  

\begin{thmx}
Suppose $1\to \Delta\to \Gamma\xrightarrow{\pi} \Lambda\to 1$ is an exact sequence of countable groups and $(N_i)$ is a sequence of finite index normal subgroups of $\Gamma$ which separates $\Delta/(\cap_i (N_i\cap \Delta))$ from $\Gamma/(\cap_i N_i)$.  This sequence induces sequences and hence box spaces of $\Delta$ and $\Lambda$.   Then if $\asdim\Box_{(N_i\cap \Delta)}\Delta = m$ and $\asdim\Box_{(\pi(N_i))}\Lambda = n$, $\asdim\Box_{(N_i)}\Gamma \leq m + n$.  
\end{thmx} 

We have therefore established a special case of the (unproved) statement \cite[Theorem A]{finnsell2015asymptotic}.  Our extra assumptions are that $\Gamma$ is countable, that the $N_i$ are normal, and that $(N_i)$ is a sequence separating $\Delta$ from $\Gamma$.  In \cite[Theorem A]{finnsell2015asymptotic}, the family of subgroups are only assumed to be semiconjugacy separating, which is weaker than assuming $(N_i)$ to be a filtration by normal subgroups.  We expect analogous techniques would work to prove these results for uncountable groups, but require translating everything to the language of coarse structures.  It is worth remarking that a countable group may have uncountably-many finite index subgroups.  

Something should be said about how our approach differs from that of \cite{finnsell2015asymptotic}.  We employ a similar general strategy in attempting to apply a Hurewicz-type theorem, and as our notation will indicate, we are using similar arguments to handle DAD as are used in \cite{Brodskiy2006AHT} to handle asymptotic dimension.  We also still need to overcome an analogous problem as in \cite[Lemma 4.1(5)]{finnsell2015asymptotic}.  An example first provided by Jiawen Zhang and restated below shows the lemma just cited is false.  We are very grateful to Jianchao Wu for helping to communicate this example and explain the problem that it highlights.  There are even `worse' counterexamples one could construct, but this one is very simple.  

\begin{theorem*}Let $6\Z$ and $9k\Z$ be subgroups of $\Z$ for $k$ odd.  Let $6\Z/(6\Z\cap 9k\Z)\to \Z/9k\Z$ be the natural inclusion (see the diagram in \cite[Remark 4.2]{finnsell2015asymptotic}).  When considered with the quotient metrics induced from a metric on $\Z$, this inclusion distorts the length of elements by an arbitrarily large amount for $k$ sufficiently large.  \end{theorem*}

\normalfont


In order to avoid this problem, we assume a kind of compatibility between the family $(N_i)$ and the particular extension.  More precisely, we assume $(N_i)$ separates $\Delta$ from $\Gamma$ (see \ref{separable subgroup definition} for a definition).  We will see in \ref{baumslag-solitar groups section} that there are still significant applications, as all elementary amenable groups are built up from extensions which are compatible with sufficiently large collections.  In fact, for our motivating example of the Baumslag-Solitar groups $\text{BS}(1, n)$, we will see that \textit{any} sequence of finite index subgroups satisfies the hypotheses of our extension theorem.  We show that box spaces of $BS(1, n)$ along filtrations have asymptotic dimension $2$ for all $n$, which in particular exhibits an infinite class of groups with exponential growth but finite-dimensional box spaces.  In obtaining the lower bound for this result, we compute the asymptotic dimension of $BS(1, n)$, which by itself involves some interesting methods in geometric group theory, specifically the notion of ends of groups, which allows one to prove that a torsion free group of asymptotic dimension $1$ is free.  As we will see, this is an immediate consequence of results already in the literature.  We also deal with some examples of more complicated groups which have similar presentations or are in some sense `generalized Baumslag-Solitar groups'.  

More generally, it is natural to try to use our extension theorem along similar lines as in \cite[Section 4]{finnsell2015asymptotic}, which attempts to relate the asymptotic dimension of box spaces of residually finite, elementary amenable groups to their Hirsch length.  While the author encountered difficulty implementing this idea because of the additional assumptions in the extension theorem, we will still discuss how the content of this paper could be used to recover similar results.  Specifically, one can probably use the Hurewicz theorem for the asymptotic dimension together with the ideas in section \ref{applications} to obtain an extension theorem for box spaces which doesn't require the collection $(N_i)$ to be a sequence; then use this to prove some special cases of the claims found in \cite[Section 4]{finnsell2015asymptotic}.

\section{(Dynamic) Asymptotic Dimension}
The following concept is originally due to Gromov, specifically \cite[Section 1.E]{Gromov1991AsymptoticIO}.  There are many definitions, but we will use one throughout which makes the analogy with the dynamic asymptotic dimension obvious:

\definition{Let $(X, d)$ be a metric space.  If $U\subset X$ is a subset, an $r$-chain in $U$ is a sequence $x_0, \ldots, x_n$ of points in $U$ such that $d(x_i, x_{i+1})\leq r$ for all $i$; and two points in $U$ are in the same $r$-component of $U$ if they are connected by an $r$-chain in $U$.  The subset $U$ is said to be $r$-connected if any two points in $U$ are connected by an $r$-chain in $U$.  Finally, $\asdim X\leq n$ means for all $r>0$ there is $M>0$ and a cover $\mathcal{U} = \{U_0, \ldots, U_n\}$ of $X$ such that each $r$-component of each $U_i$ has diameter at most $M$.  Such a cover is called a $(d, R, M)$-cover of $X$.  The asymptotic dimension of $X$ is then the smallest integer $d$ such that $\asdim X\leq d$ and is infinite if no such $d$ exists.  }

\examples{Many of the coarse spaces we consider will be graphs, considered as metric spaces by equipping their vertex set with the $l^1$-path metric.  Finitely generated groups can be considered as metric spaces through the Cayley graph construction.  More precisely, if $F\subset \Gamma$ is a finite generating set we denote by $C_F(\Gamma)$ the graph whose vertices are the elements of $\Gamma$ and which has an edge between $\gamma$ and $\delta$ if there is $f\in F$ such that $f\gamma = \delta$.  For our purposes, we can always assume $F = F^{-1}$ so that we can think of $C_F(\Gamma)$ as an undirected graph.  The $l^1$-path metric on the set of vertices is then the same as the word metric $d_F(\gamma, \delta) = |\gamma\delta^{-1}|_F$ (where $|\cdot|_F$ denotes the minimal length according to the generating set $F$).  More generally, if $\Gamma$ is a countable group, it admits a right-invariant proper metric unique up to coarse equivalence \cite[Lemma 2.3.3]{notes_on_prop_A}.  We will fix, from now on, some such metric on $\Gamma$ and note that its uniqueness up to coarse equivalence implies we can sensibly refer to $\asdim\Gamma$.  } \label{group metric}

\notation{Throughout, $\Gamma\acts X$ and $\Lambda\acts Y$ will be actions of countable groups on compact metric spaces by homeomorphisms.  }

\normalfont
We now define the dynamic asymptotic dimension of a group action, originally introduced by Guentner, Willett, and Yu in \cite{dasdimGWY}.  

\definition{Let $\Gamma\acts X$ be an action on a compact Hausdorff space by homeomorphisms.  The \textit{dynamic asymptotic dimension} of $\Gamma\acts X$, denoted $\dad(\Gamma\acts X)$, is the smallest integer $d$ such that for all finite $F\subset \Gamma$ there exists an open cover $\mathcal{U} = \{U_0, \ldots, U_d\}$ of $X$ such that, for all  $0\leq j\leq d$, the set

$$
\Bigg\{\begin{array}{l|l} g \in \Gamma & g = f_k\cdots f_1 \text{ and } \exists x\in U_j \\ & \text{ such that } f_{k_0}\cdots f_1\cdot x\in U_j \text{ }\forall \text{ } 1\leq k_0\leq k \end{array}\Bigg\}
$$

is finite.}

\normalfont
If $\Gamma\acts X$ is free, this definition is equivalent to the following: 

\definition{If $\Gamma\acts X$ is an action on a compact space by homeomorphisms, $\dad_{free}(\Gamma\acts X)$ is the smallest integer $d$ such that for all finite $F\subset \Gamma$ there is an open cover $\mathcal{U} =U_0,\ldots, U_d$ such that for all $0\leq j\leq d$ and all $x\in U_j$, the set

$$\Bigg\{\begin{array}{l|l} y\in X & \exists \text{ }  f_1, \ldots, f_k\in F \text{ such that } \\   \text{ } & y = f_k\cdots f_1\cdot x  \text{ and } f_{k_0}\cdots f_1\cdot x\in U_j \text{ } \forall \text{ } 1\leq k_0\leq k \end{array}\Bigg\}$$

is uniformly finite with respect to $x$ (i.e. the same upper bound applies for all choices of $x$).  }\label{dad definition}
\normalfont

The following terminology makes the analogy with the asymptotic dimension clear, and will be critical to the proof of our main theorem.  It is inspired by \cite{Brodskiy2006AHT}.  

\definition{Let $\mathcal{P}_{\text{fs}}(\Gamma)$ denote the collection of all finite, symmetric subsets of $\Gamma$ containing the identity.  Let $X$ be a compact metrizable space, $\Gamma\acts X$ an action by homeomorphisms, $A\subset X$ a subset, and $F\in \mathcal{P}_{\text{fs}}(\Gamma)$.  An $F$-chain in $A$ is a finite sequence $x_0, \ldots, x_n$ of points in $A$ such that $x_{i+1} = f_{i+1}\cdot x_i$ with $f_{i+1}\in F$ ($0\leq i<n$).  Two points $x, y\in A$ are in the same $F$-component of $A$ if they are connected by an $F$-chain in $A$ (this is an equivalence relation since $F$ is symmetric).  We say $A$ is $F$-connected if any two points in $A$ can be connected by an $F$-chain.  If $S\in \mathcal{P}_{\text{fs}}(\Gamma)$, we say $A\subset X$ is $S$-bounded if it is a subset of $S\cdot x$ for some $x\in X$.  A cover of $A$ by $m+1$ sets which are open (respectively, closed) in the subspace topology on $A$ and whose $F$-components are $S$-bounded is called an open (respectively, closed) $(m, F, S)$-cover of $A$.  Notice that in this case $S$ can always be taken to be $F^r$ for some $r>0$.  In case we need to distinguish between different actions on $X$ (as we will in section \ref{applications}), we will refer to an (open/closed) $(m, F, S)$-cover of $A$ with respect to $\Gamma\acts X$.  

An $n$-dimensional control function for $A$ (with respect to $\Gamma\acts X$) is a function $\mathcal{D}_{A}: \mathcal{P}_{\text{fs}}(\Gamma)\to \mathcal{P}_{\text{fs}}(\Gamma)$ such that for every $F\in \mathcal{P}_{\text{fs}}(\Gamma)$, there is an open $(n, F, \mathcal{D}_X(F))$-cover of $A$ (again, with respect to $\Gamma\acts X$ if the action is ambiguous).  Note that $\dad_{free}(\Gamma\acts X)\leq n$ iff there is an $n$-dimensional control function for $X$.  We similarly define $\dad(A)$ to be the least integer $n\geq 0$ such that there exists an $n$-dimensional control function for $A$.  } \label{chain definition}

\remark{Although the terminology in \ref{chain definition} relates specifically to the definition of $\dad$ for free actions, we will use it even for actions which are not free.  The Hurewicz-type theorem that we will prove in section \ref{product theorem} will be stated using the definition for free actions.  Note, however, that by the definition for non-free actions, any action with infinite stabilizers will be infinite dimensional.  Since the two definitions are equivalent for an action with uniformly-finite stabilizers, and actions with uniformly-finite stabilizers can easily be related to free actions by ignoring some torsion elements, we do not really lose anything by restricting our attention in this way.  From now on we will only use the definition for free actions aside from a brief mention of the general definition in section \ref{odometers and box spaces}.  } \label{freeness remark}

\normalfont
For actions on totally disconnected spaces, we can equivalently require the covers to be clopen, as shown below.  

\lemma{Suppose $\mathcal{V} = \{V_0, \ldots, V_d\}$ is an open cover of a totally disconnected, compact space $X$.  Then there is a clopen cover $\mathcal{W} = \{W_0, \ldots, W_d\}$ such that $W_i\subset V_i$ for all $0\leq i\leq d$.  In the context of the previous definition, this implies $\mathcal{W}$ is a $(d, F, S)$-cover if $\mathcal{V}$ is a $(d, F, S)$-cover.}   \label{lemma2}

\begin{proof}

Given $\mathcal{V}$, there is a cover $\mathcal{U} = \{U_i\}_{i=0}^d$ by open sets such that $\overline{U_i}\subset V_i$ for $0\leq i\leq d$ (a more general statement is shown in \cite{munkres}, Lemma 41.6).  Then since $X$ has a basis of clopen sets, we can find for each $i$ and each $x\in \overline{U}_i$ a clopen set $W_{i, x}$ containing $x$ and contained in $V_i$.  For each $i$, the $W_{i, x}$'s are a cover of $\overline{U_i}$ (which is compact since $X$ is compact); and so there is a finite subcover, the union over which we denote $W_i$.  This is a finite union of clopen sets and therefore clopen.  Moreover, the collection $\{W_i\}_{i=0}^d$ is a cover of $X$, and $W_i\subset V_i$ for all $i$.  \end{proof} 


\normalfont The following lemma appeared first in \cite[Lemma 2.9]{toprigiditydad}, and we refer there for a proof.  

\lemma{Suppose $\Gamma\acts X$ is a free action.  Then $\dad(\Gamma\acts X)\geq \asdim \Gamma$.  In fact, this inequality holds even if one modifies the definition of $\dad$ to allow covers by Borel sets, or even arbitrary sets.  }\label{easy inequality} \qed

\definition{If $\{G_n\}_{n\in \mathbb{N}}$ is a countable collection of finite metric spaces, the coarse disjoint union $\sqcup_{n} G_n$ is the metric space with underlying set the disjoint union of the $G_n$ and metric $d(x, y):=d_{G_n}(x, y)$ if $x, y\in G_n$ and $d(x, y):= \diam(G_n) + \diam(G_m)$ if $x\in G_n$ and $y\in G_m$ with $m\neq n$.  If $\Gamma$ is a group and $\{N_i\}$ is a collection of finite-index subgroups, we define $\Box_{(N_i)}\Gamma := \sqcup_i \Gamma/N_i$ where each $\Gamma/N_i$ has the metric induced from the right-invariant proper metric on $\Gamma$ we fixed earlier (i.e. the distance between cosets is the infimum over the distances between their different elements according to the original metric on $\Gamma$).  This is called a box space of $\Gamma$.  If the metric on $\Gamma$ is the word metric coming from some generating set, and if the $N_i$ are normal, the induced metrics on quotients are the word metrics coming from the image of that generating set in the quotient.  }

\definition{A filtration of $\Gamma$ is a collection $\{N_i\}$ of finite index subgroups of $\Gamma$ such that $\bigcap_i N_i = \{e\}$, the trivial subgroup.  If $\Gamma$ is residually finite, the collection of all finite index normal subgroups is then a filtration.  Note that a countable, residually finite group has a filtration which is also a sequence, that is $N_i\supset N_{i+1}$.  }\label{filtration definition}

\remark{Our definition of coarse disjoint union is slightly different from the one found in \cite[Section 2]{boxspacesDT}, and so the box space we associate to $\{N_i\}$ is not necessarily a coarse disjoint union over $\Gamma/N_i$ in their sense (as the space between different components may not tend to infinity).  However, our definition of box space is the same as theirs in the case when $\Gamma$ is finitely generated and residually finite and $\{N_i\}$ is a filtration of $\Gamma$.  Note also that we do not always require the subgroups $N_i$ to be normal, as this is not necessary for some of our results.  We will at various times need to assume that the $N_i$ are normal and/or that they form a sequence, meaning $N_{i+1}\subset N_i$.  We are therefore careful to repeat our assumptions in each statement.   }

\definition{If $(G_i)$ is a sequence of metric spaces, we say $\asdim(G_i)\leq d$ uniformly if for all $r>0$ there is an $R>0$ such that each $G_i$ admits a $(d, r, R)$-cover.  We similarly say for a sequence $(\Gamma\acts X_i)$ of actions on spaces that $\dad(\Gamma\acts X_i)\leq d$ uniformly if for all $F\in \mathcal{P}_{\text{fs}}(\Gamma)$ there is $S\in \mathcal{P}_{\text{fs}}(\Gamma)$ such that there is an open  $(d, F, S)$-cover of each $X_i$.  We can similarly define what it means for a general collection of metric spaces or actions to have $\asdim$ or $\dad$ at most $d$ uniformly.  In terms of control functions, this means having one control function which `works' for all spaces or actions in the collection.  }  \label{uniform dimension definition}

\lemma{$\asdim\Box_{(N_i)}\Gamma\leq d$ if and only if $\asdim(\Gamma/N_i)\leq d$ uniformly.  }\label{uniform and box space dimension}

\begin{proof} For one inequality, notice that a $(d, r, R)$-cover of $\Box_{(N_i)}\Gamma$ gives rise to a $(d, r, R)$-cover of each $\Gamma/N_i$ (just take the intersection). 

For the other, let $r>0$ and find a $(d, r, R)$-cover $\mathcal{U}^i = \{U^i_0, \ldots, U^i_d\}$ of each $\Gamma/N_i$.  Let $I_r = \{i\mid \diam(\Gamma/N_i)\leq r\}$.  Then $\cup_{i\in I_r} \Gamma/N_i\subset \Box_{(N_i)}\Gamma$ has diameter at most $2r$ (hence is $0$-dimensional) and the family $\{\Gamma/N_i\subset \Box_{(N_i)}\Gamma | i\notin I_r\}$ is $r$-disjoint.  The desired inequality follows from \cite[Theorem 1]{MR1808331} (we can simply put together the separate covers of the $\Gamma/N_i$ for $i\notin I_r$).  \end{proof}

\normalfont The following lemma shows that the dimension of groups and box spaces behaves nicely with respect to direct limits.  While we will not really use it, it is relevant to the discussion of elementary-amenable groups in \ref{baumslag-solitar groups section} and does not seem to be written down anywhere.  

\lemma{Suppose $G = \bigcup_i G_i$ is an increasing union of groups and $X  = \Box_{(H_j)} G = \bigcup_i X_i$ is an increasing union of box spaces where $X_i = \Box_{(H_j\cap G_i)}G_i$ is a box space of $G_i$ (here the metrics on $X$ and $G_i$ are induced from the right invariant, proper metric on $G$, and the metric on $X_i$ is constructed using the metric on $G_i$).  Then $\asdim G = \sup_i \asdim G_i$ and $\asdim X = \sup_i \asdim X_i$.  }\label{direct limit lemma}

\begin{proof}
Fix $R>0$.  Then there is $N$ such that $B_e^R(G)\subset G_N$.  Notice that $d_G(G_N\cdot g, G_N\cdot h)$ is either $0$ or $>R$.  This follows since if $d_G(G_N\cdot g, G_N\cdot h)\leq R$ then the two orbits have a nontrivial intersection and therefore coincide.  Since different $G_N$-orbits are $>R$ apart, we can cover each one separately using $\asdim G_N + 1$ sets.  For the second claim, find $N$ as before.  Then if $x, y\in G/H_j$ are in separate $G_N$-orbits, their representatives in $G$ must be in separate $G_N$ orbits, and so $d_{G/{H_j}}(x, y)>R$.  We can therefore similarly cover each orbit separately using $\asdim \Box_{(H_j\cap G_N)} G_N + 1$ sets.  \end{proof}

\normalfont
In addition to \ref{easy inequality}, we see that when the space being acted on is discrete, there is a very concrete relationship between the dynamic asymptotic dimension of the action and the asymptotic dimension of the orbits.  

\lemma{Suppose $\Gamma$ is a countable group with subgroup $N<\Gamma$ and recall we have equipped $\Gamma$ with a right-invariant, proper metric (see \ref{group metric}).  Suppose further that $\Gamma/N$ is a quotient of $\Gamma$ with the metric induced from $\Gamma$.  Consider the action $\Gamma\acts \Gamma/N$ on left cosets by left multiplication (this is an action even if $N$ is not normal).  Then an $S$-bounded subset has diameter at most $\diam(S)$ in $\Gamma$.  Moreover, a subset $A\subset \Gamma/N$ with diameter at most $R$ is $B_e^R(\Gamma)$-bounded (recall that quotients do not increase the growth of a group).  In particular, this implies a $(d, B_e^r(\Gamma), S)$-cover of $\Gamma/N$ is a $(d, r, \diam(S))$-cover of $\Gamma/N$; and a $(d, r, R)$-cover of $\Gamma/N$ is a $(d, B_e^r(\Gamma), B_e^R(\Gamma))$-cover of $\Gamma/N$.  }\label{basic lemma} 

\begin{proof}
Observe that $x, y\in A\subset \Gamma$ are connected by a $B_e^r(\Gamma)$-chain in A iff they are connected by an $r$-chain in $A\subset \Gamma$.  This follows since $d_{\Gamma/N}(\gamma\cdot x, x)\leq r$ for any $\gamma\in B_e^{r}(\Gamma)$ and conversely any elements $\leq r$ away from each other in $\Gamma/N$ are connected by some $\gamma\in B_e^r(\Gamma)$.  \end{proof}

\normalfont
There is a familiar dynamical object which is reminiscent of a box space called an odometer.  We will see in \ref{odometers and box spaces} that the asymptotic dimension of box spaces is closely related to the dynamic asymptotic dimension of odometers.  

\definition{Let $\Gamma$ be a countable group and $(N_i)$ ($i\in I$) a countably infinite collection of finite-index subgroups directed by inclusion.  The profinite completion of $\Gamma$ along $(N_i)$ is the inverse limit $\displaystyle{\lim_{\leftarrow i} \Gamma/N_i}$ and we denote it $\widehat{\Gamma}_{(N_i)}$.  If $(N_i)$ is the collection of all such subgroups, we just write $\widehat{\Gamma}$.  If the $N_i$ are normal subgroups, there is a group structure on $\widehat{\Gamma}_{(N_i)}$ and an action $\Gamma\acts \widehat{\Gamma}_{(N_i)}$ given by the quotient map and left multiplication.  Otherwise, there is still an action $\Gamma\acts \widehat{\Gamma}_{(N_i)}$ by left multiplication on left cosets in each coordinate.  The canonical quotient map $p_i: \widehat{\Gamma}_{(N_i)}\to \Gamma/N_i$ is continuous for the profinite topology and equivariant.  The profinite topology can also be thought of as the subspace topology inherited from the Tychonoff topology on $\prod_i \Gamma/N_i$.  Notice that $\widehat{\Gamma}_{(N_i)}$ is closed in $\prod_i \Gamma/N_i$.  In fact, if $\phi: \Gamma\to \prod_i \Gamma/N_i$ is the map coming from all the quotient maps, $\widehat{\Gamma}_{(N_i)} = \overline{\phi(\Gamma)}$ (using that $(N_i)$ is directed).  }


\section{Odometers and Box Spaces} \label{odometers and box spaces}
\normalfont
As we have seen, asymptotic dimension and DAD are closely related at least when the space being acted on is discrete.  This relationship is encapsulated by \ref{basic lemma} and, together with some topological facts about profinite completions, can be used to show the DAD (according to the definition for free actions) of an odometer formed from a sequence $(N_i)$ of subgroups is the same as the asymptotic dimension of the box space formed from $(N_i)$.  Consequently, we see that questions about asymptotic dimension of such box spaces are special cases of questions about $\dad$.  This connection is spiritually similar to how the DAD of an action is bounded above by the asymptotic dimension of its warped cone (this can be shown directly, but see \cite{warpedcones} for definitions and a more complete discussion of warped cones and DAD).  

%
%
%

\theorem{Let $\Gamma$ be a countable group and $(N_i)$ ($i\in I$) a countable collection of finite-index subgroups of $\Gamma$ directed by inclusion.  Then $\dad_{free}(\Gamma\acts \Gamma_{(N_i)}) \leq \asdim \Box_{(N_i)} \Gamma$.  } \label{dimension of odometers and box spaces}

\begin{proof}


Suppose $\asdim\Box_{(N_i)}\Gamma\leq d$ and fix a finite subset $F\in \Gamma$.  Find $R>0$ so that $d_\Gamma(\gamma\cdot x, x) < R$ for all $\gamma\in F$ and $x\in \Gamma$, hence $d_{\Gamma/N_i}(\gamma\cdot x, x)< R$ for all $\gamma\in F$ and all $i$.  Suppose $\asdim\Box_{(N_i)}\Gamma \leq d$.  Then $\asdim(\Gamma/N_i)\leq d$ uniformly.  Let $\mathcal{U}^i = \{U^i_0, \ldots, U^i_n\}$ be a $(d, R, M_R)$-cover of $\Gamma/{N_i}$.  Then by \ref{basic lemma} each $\mathcal{U}^i$ is also a $(d, F, B_{e}^{M_R + 1}(\Gamma/N_i))$-cover of $\Gamma/{N_i}$ (and since the cardinality of balls in finite quotients of $\Gamma$ is bounded by the cardinality of balls in $\Gamma$, this bound is independent of $i$).  If $K = \ker(\Gamma\to \widehat{\Gamma}_{(N_i)})$, then $\Gamma/K\acts \widehat{\Gamma}_{(N_i)}$ freely and is orbit equivalent to the original $\Gamma$-odometer.  Find $i_F$ so that the quotient $\Gamma/K\to \Gamma/N_{i_F}$ induces a map $C_F(\Gamma/K)\to C_F(\Gamma/N_{i_F})$ which is an isometry on balls of radius $\leq M_R + 2$ (this is possible, for isntance, by \cite[Proposition 2.1]{boxspacesDT}).  Then pulling back $\mathcal{U}^{i_F}$ by the projection $p_{i_F}: \widehat{\Gamma}_{(N_i)}\to \Gamma/N_{i_F}$ gives a $(d, F, B_e^{M_R + 1}(\Gamma/N_{i_F}))$-cover of $ \widehat{\Gamma}_{(N_i)}$.  \end{proof}


\normalfont

It will be important for \ref{prop2} that the collection $(N_i)$ is a sequence, meaning the index set is the natural numbers and $N_{i+1}\subset N_i$.  


\lemma{Let $\Gamma$ be a countable group and $(N_i)$ a countable collection of finite index subgroups which is directed by inclusion.  A subset $U\subset \widehat{\Gamma}_{(N_i)}$ is clopen iff it is of the form $p_i^{-1}(S)$ where $S\subset \Gamma/N_i$ is any subset. } \label{lemma3}

\begin{proof}

First, sets of the form $p_i^{-1}(S)$ are clopen as $S\subset \Gamma/N_i$ is clopen and $p_i$ is continuous.  

Fix the metric $d((x_i), (y_i)) = \sum_{i=1}^{\infty} \frac{d_{\Gamma/N_i}(x_i, y_i)}{2^i\diam(\Gamma/N_i)}$ on $\widehat{\Gamma}$.  If $U$ is clopen, then $\widehat{\Gamma} = U\cup U^c$ where $U$ and $U^c$ are separate connected components.  In fact, since $\widehat{\Gamma}$ is a compact metric space, there is $\epsilon$ such that the $\epsilon$-neighborhoods of $U$ and $U^c$ do not intersect.  

Now, suppose $U$ is not of the desired form.  Then for all $i$ there is $j\geq i$ and $x_j\in\widehat{\Gamma}$ such that $p_j(x_j)\in p_j(U)$ but $p_{j+1}(x_j)\notin p_{j+1}(U)$.  But then $x_j\notin U$ for all such (inifinitely-many) $j$ and $d(x_j, U)\leq \sum_{n=j+1}^\infty \frac{1}{2^n}\to 0$ as $j\to \infty$\footnote{Technically, the sum is over a rearrangement, but that doesn't matter.}.  Hence, there are points in $U^c$ arbitrarily close to $U$.  By what was said above, this implies $U$ is not clopen. \end{proof}

\proposition{Let $\Gamma$ be a countable group and $(N_i)$ a sequence of finite-index subgroups which is directed by inclusion (so $i\in \mathbb{N}$ and $N_{i+1}\subset N_i$).  Let $\Gamma\acts \widehat{\Gamma}_{(N_i)}$ be the odometer action induced by left multiplication on each finite quotient.  If $\dad_{free}(\Gamma\acts \widehat{\Gamma}_{(N_i)})\leq d$ then $\asdim(\Box_{(N_i)}\Gamma)\leq d$.  } \label{prop2}

\begin{proof}
Fix a proper, left-invariant metric on $\Gamma$.  Let $R\geq 0$ be given and set $F = B_e^R(\Gamma)$.  Let $\mathcal{V} = \{V_0, \ldots, V_d\}$ be an open $(d, F, S)$-cover of $\widehat{\Gamma}_{(N_i)}$.  Observe that for each $i$, the distance between the $F$-components of $p_i(V_j)$ is $>R$ according to the metric on $\Gamma/N_i$.  



Now by \ref{lemma2} we can assume each $V_j$ is clopen and hence by \ref{lemma3} is the pullback of $T\subset \Gamma/N_i$ for some $i$.  Find $i_0$ such that all of the $V_j$ are pullbacks of subsets of $\Gamma/N_{i_0}$.  Then if $i\geq i_0$, $p_i(x)\in p_i(V_j)$, and $f\cdot p_i(x)\in p_i(V_j)$, then $x, f\cdot x\in V_j$.  Thus, each $F$-component of $p_{i}(V_j)$ is $S$-bounded for all $i\geq i_0$; so by \ref{basic lemma} the diameter (as a subset of $\Box_{(N_i)}\Gamma$) of each $R$-component of $p_i(V_j)\subset \Gamma/N_i$ is bounded by $\diam(S)$ for $i\geq i_0$.  This gives a $(d, R, \diam(S))$-cover of each $\Gamma/N_i$ for all $i\geq i_0$.  Since $(N_i)$ is a sequence, the set $\{i \mid i < i_0\}$ is finite, and so the diameter of $\Gamma/N_i$ for $i<i_0$ is uniformly bounded.  We have therefore shown $\asdim(\Gamma/N_i)\leq d$ uniformly.  The last part then follows from \ref{uniform and box space dimension}.   \end{proof}

\normalfont
We now state both inequalities (when $(N_i)$ is a sequence) in one theorem:

\theorem{Let $\Gamma$ be a countable group and $(N_i)$ a sequence of finite index subgroups which is directed by inclusion.  Suppose $\Gamma\acts \widehat{\Gamma}_{(N_i)}$ by left multiplication.  Then $\dad_{free}(\Gamma\acts \widehat{\Gamma}_{(N_i)}) = \asdim\Box_{(N_i)}\Gamma$.  } \label{box space theorem} \qed


\corollary{Let $\Gamma$ be a countable, residually finite group and $(N_i)$ a sequence of finite index subgroups which is also a filtration of $\Gamma$ (see \ref{filtration definition}).  Then $\dad(\Gamma\acts \widehat{\Gamma}_{(N_i)}) = \asdim\Box_{(N_i)}\Gamma$.  } \qed


\normalfont
This can also be combined with \cite[Proposition 3.1]{boxspacesDT}:

\corollary{Let $\Gamma$ be a finitely generated, residually finite group and $(N_i)$ a sequence of finite index subgroups which is also a filtration of $\Gamma$.  Then either $\dad(\Gamma\acts \widehat{\Gamma}_{(N_i)}) = \infty$ or $\dad(\Gamma\acts \widehat{\Gamma}_{(N_i)}) = \asdim \Gamma$.  The dimension is finite if $\Gamma$ has polynomial growth (equivalently is virtually nilpotent).  } \qed

%

\section{Product Theorem}\label{product theorem}
\normalfont
Now we begin our study of permanence properties of $\dad$ which will lead to our main result.  Most basically, one would expect DAD to be subadditive over products of groups acting on products of spaces, and the proof of this fact will exhibit some of the same ideas used in the next, more technical, section.  A product theorem for the dynamic asymptotic dimension of groupoids was recently established by B\"{o}nicke in \cite{bonicke2023dynamic}, but we will go on to prove the more general Hurewicz mapping theorem for the transformation group case.  It is worth remarking that the product theorem for groupoids does not really require any additional ideas from the case for group actions.  

We rely heavily on the ideas of \cite{Brodskiy2006AHT} for both this section and the next and have noted when one of our lemmas is inspired by one of theirs.  Indeed, if one relaxes the definition of $\dad$ to allow Borel sets, or even sets obtained from finite operations with open sets, the work of \cite{Brodskiy2006AHT} can be transported to the dynamical setting with very little extra work.  Our solution to the extra technicality created by the topology of $X$ comes from a lemma provided in this section which we refer to as `the trick'.  This lemma capitalizes on some rigidity to provide a way of switching between covers using open or closed sets while preserving the properties related to $\dad$.  Part (2) was previously established in \cite[Lemma 2.8(i)]{toprigiditydad}, and the essential content of part (1) also appeared in the surrounding discussion in that paper.    

\definition{Let $\Gamma\acts X$ be an action on a compact metric space by homeomorphisms and $A\subseteq X$.  If $k\geq n+1\geq 1$, an $(n, k)$-dimensional control function for $A$ (with respect to $\Gamma\acts X$) is a function $\mathcal{D}_{A}:\mathcal{P}_{\text{fs}}(\Gamma)\to \mathcal{P}_{\text{fs}}(\Gamma)$ such that for any $F\in \mathcal{P}_{\text{fs}}(\Gamma)$ there is a $(k-1, F, \mathcal{D}_{A}(F))$-cover $\mathcal{U} = \{U_0, \ldots, U_{k-1}\}$ of $A$ (with respect to $\Gamma\acts X$) by sets open in the subspace topology on $A$ such that each $x\in A$ belongs to at least $k-n$ elements of $\mathcal{U}$ (equivalently $\cup_{i\in T} U_i = A$ for every $T\subset \{0, \ldots, k-1\}$ consisting of $n+1$ elements). } \label{control function space definition}


\lemma{(the trick) Let $\Gamma\acts X$ be an action on a compact metric space by homeomorphisms and $A\subset X$ closed.  Let $F, S\in\mathcal{P}_{\text{fs}}(\Gamma)$.  
\begin{itemize}
\item[(1)] If $n\geq0$, $k\geq n+1$, and $\mathcal{U} = \{U_0, \ldots, U_{k-1}\}$ is a $(k-1, F, S)$-cover of $A$ by sets open in the subspace topology on $A$ such that each $x$ belongs to $\geq k-n$ elements.  Then there is a $(k-1, F, S)$-cover by sets which are the closures of open sets (so closed in $X$) and still have the property that each $x$ belongs to $\geq k-n$ elements.  

\item[(2)] If the $F$-components of $A$ are $S$-bounded, then $A$ is contained in a set which is open in $X$ and whose $F$-components are $S$-bounded.  

\end{itemize} } \label{control functions lemma}

\begin{proof}
For the first part, assume $k>1$ and let $\mathcal{U} = \{U_0, \ldots, U_{k-1}\}$ be such a cover.  Then $\bigcup_{S\subset \{0, \ldots, k\} : |S| = k-n} \cap_{s\in S}U_s$ is an open cover, and so has some Lebesgue number $\lambda$.  This shows every $x\in X$ has a $\lambda$-ball which is contained in at least $k-n$ elements of $\mathcal{U}$.  Replace the $U_i$ with the closures of their $\lambda/2$-interiors.  If $n=0$ and $k=1$, $\mathcal{U} = \{A\}$, which is already a cover by closed sets (closed in $X$ in fact).  

%
%
%
To prove $(2)$, suppose not, that is, suppose for every $n$ that the $F$ components of $N_{1/n}(A)$ are not $S$-bounded.  So for every $n$, we have an $F$-chain $x^{(n)}_0, \ldots, x^{(n)}_{K_n}$ in $N_{1/n}(A)$ with $x^{(n)}_0, \ldots, x^{(n)}_{K_n-1}\in S\cdot x_0$ but $x_{K_n} = \gamma_n\cdot x^{(n)}_0$ for $\gamma_n = f_ns_n$ with $s_n\in S$, $f_n\in F$, and $\gamma_n\notin S$.  Since there are finitely-many $\gamma\in \Gamma$ with the property that $\gamma = fs$ for $s\in S$ and $f\in F$, we must have $\gamma_n = \lambda$ for infinitely-many $n$.  Furthermore, there are finitely-many sequences $f_1, \ldots, f_{K_n}$ ($f_k\in F$) such that $f_k\cdots f_1\in S$ for all $1\leq k<K_n$ but $f_{K_n}\cdots f_1 = \lambda\notin S$, and so for infinitely-many $n$ we have an $F$-chain $x^{(n)}_0, \ldots, x^{(n)}_{K}$ in $N_{1/n}(A)$ and $f_1, \ldots, f_K\in F$ such that $f_{k+1}\cdot x^{(n)}_{k} = x^{(n)}_{k+1}$ for $0\leq k<K$, and also $x_K = \lambda\cdot x_0$ with $\lambda\notin S$.  Since $X$ is compact, we have a subsequence $x_0^{(n_j)}$ converging to a point $y_0\in X$, and so $f_k\ldots f_1\cdot x_0^{(n_j)} = x_k^{(n_j)}$ converges to $f_k\ldots f_1\cdot y_0$ for all $1\leq k\leq K$ by continuity.  Since $x_k^{(n_j)}\in N_{1/n_j}(A)$ and $A$ is closed in $X$, $x_k^{(n_j)}$ converges to a point in $A$.  Hence, we have an $F$-chain in $A$ which is not $S$-bounded.  \end{proof}

\remark{We will frequently apply the lemma above by simply saying `use the trick', when we want to replace an open cover by a closed one or a closed set with a larger open set.  For a closed subspace $A\subset X$, notice that the trick implies it is equivalent to define $\dad(A)$ and related terms using covers by sets which are open in $X$ (the definition in \ref{chain definition}, used sets which are open in the subspace topology on $A$).} \label{trick remark}

\lemma{(compare to \cite[Theorem 2.4]{Brodskiy2006AHT}) Suppose $A\subset X$ is closed and that $\mathcal{D}_{A}^{(n+1)}$ is an $n$-dimensional control function for $A$.  Define $\{\mathcal{D}_{A}^{(i)} \}_{i\geq n+1}$ inductively on powers of $F$ by $\mathcal{D}_{A}^{(i+1)}(F^r) := F^r\cdot\mathcal{D}^{(i)}_{A}(F^{3r})\cdot F^r$; then for each $k$, there is an open $(k-1, F^r, \mathcal{D}^{(k)}_{A}(F^r))$-cover of $A$ such that each $x\in A$ is covered by at least $k - n$ elements of the cover.  (We say it this way since the functions $\mathcal{D}_{A}^{(k)}$ are not defined on all elements of $\mathcal{P}_{\text{fs}}(\Gamma)$ unless $\Gamma$ is finitely generated, so we can't quite call them control functions).  }\label{extra cover lemma}

\begin{proof}

The case $k=n+1$ is obvious.  Suppose the result holds for some larger $k$.  Let $\mathcal{U} = \{U_1, \ldots, U_k\}$ be an open $(k-1, F^{3r}, \mathcal{D}_{A}^{(k)}(F^{3r}))$-cover of $A$ with all $x\in A$ in at least $k-n$ elements of $\mathcal{U}$.  Using (1) from the trick (see \ref{trick remark}), we can replace the sets $U_i$ by the closures of open sets, $V_i$, which form a cover of $A$ with the same properties.  Define $U_i' = (F^r\cdot U_i)\cap A$ (these are closed sets).  Then $\mathcal{U}' = \{U_1', \ldots, U_k'\}$ is a closed $(k-1, F^r, F^r\cdot \mathcal{D}_{A}^{(k)}(F^{3r})\cdot F^r)$-cover of $A$.  Recall that $F^r\cdot \mathcal{D}_{A}^{(k)}(F^{3r})\cdot F^r = \mathcal{D}_{A}^{(k+1)}(F^{r})$ by definition.  


Define $U_{k+1}'$ to be the union of all the (disjoint, open) sets $W_S := (\cap_{s\in S} V_s)\setminus \cup_{i\notin S} U_i'$ where $S\subset \{1, \ldots, k\}$ has $k-n$ elements.  We claim that for $S\neq T$ (where $S, T\subset \{1, \ldots, k\}$ each have $k-n$ elements), the sets $W_S$ and $W_T$ are $F^r$-disjoint (i.e. not connected by an $F^r$-chain contained in their union).  

Suppose $a, b\in U_{k+1}'\subset A$ are such that $a\in \cap_{t\in T} V_t \setminus \cup_{i\notin T} U_i'$ and $b\in \cap_{s\in S} V_s \setminus \cup_{i\notin S} U_i'$ for $T\neq S$ and $|S| = |T| = k-n$.  Then there is $t\in T\setminus S$ such that $a\in V_t$, so if $f\cdot a = b$ for some $f\in F^r$, then $b\in U_t'\supseteq F^r\cdot V_t$, a contradiction.  

Now suppose $x\in A$ belongs to exactly $k-n$ sets $U_i'$, $i\leq k$.  Since the sets $\cap_{s\in S}V_s$ form a cover as $S$ ranges over subsets of $\{1, \ldots, k\}$ with $|S| = k-n$; $x\in \cap_{s\in S}V_s$ for some $s\in S$.  But then if $x\notin U_{k+1}'$, we must have $x\in U'_j$ for some $j\notin S$, contradicting that $x$ is only contained in $k-n$ sets $U_i'$.  

Thus, $\{U'_{k+1}, U'_k, \ldots, U'_1\}$ is a $(k, F^r, \mathcal{D}_{A}^{(k+1)}(F^{r}))$-cover with each $x\in A$ contained in at least $k-n + 1$ sets (recall that $U_{k+1}'$ is open and the other sets are closed).  We can then apply the trick (2) to each closed set to produce a cover by open sets with the same properties.  \end{proof} 


\theorem{(compare to \cite[Theorem 2.5]{Brodskiy2006AHT}) If $\Gamma\acts X$ and $\Lambda\acts Y$ are free, then $\text{DAD}(\Gamma\times \Lambda\acts X\times Y)\leq \text{DAD}(\Gamma\acts X) + \text{DAD}(\Lambda\acts Y)$.  

\begin{proof}

Fix $F\in\mathcal{P}_{\text{fs}}(\Gamma)$ and $H\in\mathcal{P}_{\text{fs}}(\Lambda)$.  Suppose $\Gamma\acts X$ and $\Lambda\acts Y$ have dimension at most $m$ and at most $n$, respectively, and let $\mathcal{D}_X$ and $\mathcal{D}_Y$ be respective control functions.  Let $k = m + n + 1$.  Use \ref{extra cover lemma} to produce the functions $\mathcal{D}^{(k)}_{X}$ and $\mathcal{D}^{(k)}_{Y}$.  So there is an open $(k-1, F, \mathcal{D}^{(k)}_{\Gamma\acts X}(F))$-cover $\mathcal{U} = \{U_i\}_{i=0}^{k-1}$ of $X$ and an open $(k-1, H, \mathcal{D}^{(k)}_{\Lambda\acts Y}(H))$-cover $\mathcal{V} = \{V_i\}_{i=0}^{k-1}$ of $Y$ such that each $x\in X$ is in at least $k-m$ elements of $\mathcal{U}$ and each $y\in Y$ is in at least $k-n$ elements of $\mathcal{V}$.  Since $k-m + k - n = k+1$, the family $\{U_i\times V_i\}_{i=0}^{k-1}$ is a cover of $X\times Y$.  Moreover, it is a $(F\times H, \mathcal{D}^{(k)}_{\Gamma\acts X}(F) \cdot \mathcal{D}^{(k)}_{\Lambda\acts Y}(H))$-cover of $\Gamma\times \Lambda\acts X\times Y$.  \end{proof}

\section{Hurewicz-type Theorem}
\normalfont
In this section, we will prove a more general permanence property witnessing a natural dimension-theoretic property of $\dad$.  We start by defining the dimension of an appropriate morphism in the category of group actions.  As dynamic asymptotic dimension is meant to measure the large scale complexity of orbits, the appropriate kind of morphism here is one which sends orbits to orbits and induces a coarse map of the orbit structures.  The following definition is not new, but it was inspired by the definition of continuous orbit couple in \cite[Definition 2.6]{LiDynamics}.  A continuous orbit couple induces a `coarse equivalence of orbits', so we weaken this definition to one which will give a `coarse map' of orbits.  

\definition{Suppose $\Gamma\acts X$ and $\Lambda\acts Y$ are actions on compact metric spaces by homeomorphisms and that $f: X\to Y$ and $\alpha: \Gamma \times X\to \Lambda$ are continuous maps with $f(\gamma\cdot x) = \alpha(\gamma, x)\cdot x$.  Define $\dad(f) = \sup_{A\subset Y : \dad(A) = 0} \dad_{free}(f^{-1}(A))$ (with $\dad(f^{-1}(A))$ defined as in \ref{chain definition}).  }\label{dad of a map}

\normalfont
The intuition behind this definition is we want to measure the largest dimension of a subspace of $X$ which $f$ collapses down to something zero dimensional in the image.  Recall for instance that in the category of vector spaces, this would be the nullity of a linear map, i.e. the dimension of the pullback of $\{0\}$.  In the category of topological spaces, one can similarly define the dimension of a function to be the supremum over the dimensions of pullbacks of points.  However, as in the case of asymptotic dimension, we must alter this definition slightly by considering pullbacks of arbitrary $0$-dimensional subspaces.  See  \cite[Proposition 4.1]{Brodskiy2006AHT} for an example which shows this deviation is necessary.  



\normalfont
Many of the ideas in this section come from \cite{Brodskiy2006AHT}, so as in the previous section we note at the beginning of each statement what lemma is being analogized.  

\definition{Let $\Gamma\acts X$, $\Lambda\acts Y$ and $f: X\to Y$ be as in \ref{dad of a map}.  An $m$-dimensional control function for $f$ is a function $\mathcal{D}_f: \mathcal{P}_{\text{fs}}(\Gamma)\times \mathcal{P}_{\text{fs}}(\Lambda)\to \mathcal{P}_{\text{fs}}(\Gamma)$ such that for all $S, T\in \mathcal{P}_{\text{fs}}(\Gamma)$, we have that if $B\subset Y$ is $T$-bounded, then $f^{-1}(B)$ has a $(m, S, \mathcal{D}_f(S, T))$-cover (by sets which are open in the subspace topology on $f^{-1}(B)$, see \ref{trick remark}).  }

\lemma{(see \cite[Proposition 4.5]{Brodskiy2006AHT}) Let $\Gamma\acts X$, $\Lambda\acts Y$ and $f: X\to Y$ be as in \ref{dad of a map}.  If $\dad(f)\leq m$, then there is an $m$-dimensional control function for $f$.  }\label{dimension implies function}

\begin{proof}
Fix $T\in \mathcal{P}_{\text{fs}}(\Lambda)$ and $S\in \mathcal{P}_{\text{fs}}(\Gamma)$.  Suppose for a contradiction that $F_n$ is an increasing sequence of sets in $\mathcal{P}_{\text{fs}}(\Gamma)$ with $\cup_n F_n = \Gamma$, and that for each $n$, $B_n = T\cdot x_n$ is a $T$-bounded set such that $f^{-1}(B_n)$ does not have an open $(m, S, F_n)$-cover.  Then for any subsequence $n_k$, $f^{-1}(\bigcup_{k} B_{n_k})$ must have $\dad$ greater than $m$.  

Notice that at least one of the following alternatives must hold: either there exist infinitely many $n$ with $x_n$ in different orbits, or there exist infinitely many $n$ with $x_n$ in the same orbit.  In the first case, we get a subsequence $n_k$ with each $x_{n_k}$ in a different orbit so that $\bigcup_{k} B_{n_k}$ is zero-dimensional; this is a contradiction.  In the second case, we have a subsequence $n_k$ with each $x_{n_k}$ in the same orbit, say the orbit of $x_0$.  Fix a right-invariant, proper metric $d_\Lambda$ on $\Lambda$.  Then define a proper metric $d$ on the orbit of $x_0$ by the formula $d(x, y):=\inf_{\{\lambda\in \Lambda : \lambda\cdot x = y\}} d_\Lambda(\lambda, e)$.  Then either $d(x_{n_k}, x_0)\to \infty$ or these distances do not tend to infinity.  In the case of the latter, that implies $\bigcup_{k} B_{n_k}$ is finite and hence zero dimensional, a contradiction.  In the case of the former, passing to a further subsequence yields a sequence $n_l$ with $d(x_{n_l}, \bigcup_{j<l} x_{n_j})\to \infty$.  Again, this implies $\bigcup_{l} B_{n_l}$ is zero dimensional, which is a contradiction.  \end{proof}


\lemma{Suppose $X$ and $Y$ are compact metric spaces and $f: X\to Y$ is continuous.  Let $C\subset Y$ be closed and $N$ an open neighborhood of $f^{-1}(C)$.  Then there is an open neighborhood $G$ of $C$ such that $f^{-1}(G)\subset N$.  } \label{neighborhoods lemma}

\begin{proof}
Using the metric on $Y$, one can construct a family $\{G_\alpha\}$ of open sets containing $C$ such that $\bigcap_\alpha \overline{G}_\alpha = C$ ($\overline{G}_\alpha$ denotes the closure of $G_\alpha$).  Then $\bigcap_\alpha f^{-1}(\overline{G}_\alpha) = f^{-1}(\bigcap_\alpha \overline{G}_\alpha) = f^{-1}(C)$.  That is, for all $x\in f^{-1}(C)^c$, there is $\alpha$ such that $x\notin f^{-1}(\overline{G}_\alpha)$.  Moreover, since $f^{-1}(\overline{G}_\alpha)$ is closed, $d(x, G_{\alpha})>\delta_x$ for some $\delta_x>0$ 

Now let $N\supset f^{-1}(C)$ be an open neighborhood.  Then $N^c$ is closed and therefore compact.  For each $x\in N^c\subset f^{-1}(C)^c$, let $\delta_x>0$ as above.  Then the collection of open balls $\{B_x^{\delta_x}(X)\}$ is a cover of $N^c$ and therefore has a finite subcover $\{B_{x_k}^{\delta_{x_k}}(X)\}_{k=1}^K$.  Let $G_{\alpha_k}$ be the neighborhoods from the paragraph above coming from $x_k$.  Then if we put $G = \bigcap_k G_{\alpha_k}$, $G$ is an open neighborhood of $C$ and $f^{-1}(\bigcap_{\alpha_k}G_{\alpha_k}) = \bigcap_{\alpha_k} f^{-1}(G_{\alpha_k})\subset N$.  \end{proof}

%
%

\remark{To prove the main theorem in this section, it is enough for \ref{neighborhoods lemma} to hold for $C$ a finite subset.  }

\lemma{(compare to \cite[Proposition 4.7]{Brodskiy2006AHT}) Let $\Gamma\acts X$, $\Lambda\acts Y$, and $f: X\to Y$ be as in \ref{dad of a map}.  Suppose $\mathcal{D}_f: \mathcal{P}_{\text{fs}}(\Gamma)\times \mathcal{P}_{\text{fs}}(\Lambda)\to \mathcal{P}_{\text{fs}}(\Gamma)$ is an $m$-dimensional control function for $f$.  If we define $\{\mathcal{D}_f^{(i)}\}_{i\geq m+1}$ inductively on powers of $F\in \mathcal{P}_{\text{fs}}(\Gamma)$ by $\mathcal{D}_f^{(i+1)}(F^r, T):= F^r\cdot \mathcal{D}_f^{(i)}(F^{3r}, T)\cdot F^r$, then for each $k$ and $T$-bounded subset $B\subset Y$, there is an open $(k-1, F^r, \mathcal{D}_f^{(k)}(F^r, T))$-cover of $f^{-1}(B)$ such that each $x\in X$ is covered by at least $k-m$ elements of the cover.  }\label{function extra cover lemma}

\begin{proof}
The function $F\mapsto \mathcal{D}_f(F, T)$ is an $m$-dimensional control function for $f^{-1}(B)$; more precisely, it is a control function for \textit{any} $T$-bounded subset $B$ (note that such a subset is necessarily closed).  We can therefore apply \ref{extra cover lemma} to produce $(m, k)$-dimensional control functions for $f^{-1}(B)$ (technically these functions are only defined on powers of $F$, but that's all we need to establish this lemma).  Since the construction in the proof of \ref{extra cover lemma} depends only on the control function we start with, these new control functions can be taken to be the same for any $T$-bounded subset $B\subset Y$, and so we have the desired conclusion.  \end{proof}

\definition{Let $f: X\to Y$ be as in \ref{dad of a map}.  Suppose that $F\in \mathcal{P}_{\text{fs}}(\Gamma)$ and $S\in \mathcal{P}_{\text{fs}}(\Lambda)$.   An $(F, S)$-chain in $U$ is an $F$-chain whose image under $f$ is an $S$-chain.  Two points $x_1, x_2\in U$ are in the same $(F, S)$-component if they are connected by an $(F, S)$-chain.  We say a subset $A\subset X$ is $(F, S)$-bounded if $A$ is $F$-bounded and $f(A)$ is $S$-bounded.  } \label{more chain definition}

\lemma{(see \cite[Proposition 4.8]{Brodskiy2006AHT}) Let $\Gamma\acts X$, $\Lambda\acts Y$ and $f: X\to Y$ be as in \ref{dad of a map}.  Suppose $\mathcal{D}_f$ is an $m$-dimensional control function for $f$.  Then for all $T, S\in \mathcal{P}_{\text{fs}}(\Lambda)$, $F\in \mathcal{P}_{\text{fs}}(\Gamma)$, and $B\subset Y$ an open set with $T$-components $S$-bounded, $f^{-1}(B)$ can be covered by $k$ sets which are open in $X$ and whose $(F, T)$-components are $\mathcal{D}^{(k)}_{f}(F, S)$-bounded (with $\mathcal{D}_f^{(k)}$ defined as in \ref{function extra cover lemma}), and every $x\in f^{-1}(B)$ is contained in at least $k-m$ elements of this cover.  } \label{control function tricky lemma}

\begin{proof}
Let $C_{\alpha}$ be a $T$-component of $B$.  By assumption, $C_{\alpha} = S_\alpha\cdot x_\alpha$ where $S_\alpha\subset S$, so by \ref{function extra cover lemma}, $f^{-1}(C_\alpha)$ has a $(k-1, F, \mathcal{D}^{(k)}_f(F, S))$-cover by sets which are open in the subspace topology on $f^{-1}(C_\alpha)$ and which has the additional property that every element of $f^{-1}(C_\alpha)$ is contained in at least $k-m$ elements of the cover.  Since $f^{-1}(C_\alpha)$ is closed and hence compact, we can use the trick (1) to replace this cover with a cover by closed sets (which are closed in $X$), and then use the trick (2) to again replace the cover with a cover by sets which are open in $X$, all the while preserving the fact that every element of $f^{-1}(C_\alpha)$ is contained in at least $k-m$ elements of this new cover.  Call this cover $\mathcal{U}^\alpha = \{U^\alpha_0, \ldots, U^\alpha_{k-1}\}$.  It covers not just $f^{-1}(C_\alpha)$, but an open neighborhood of it which we can assume without loss of generality to be contained in $f^{-1}(B)$.  By \ref{neighborhoods lemma}, it is possible to find an open neighborhood $G_\alpha$ of $C_\alpha$ that $f^{-1}(G_\alpha)$ is covered by $\mathcal{U}^\alpha$.  



Find an open neighborhood $M_\alpha$ of $x_\alpha$ such that $N_\alpha := S_\alpha\cdot M_\alpha\subset G_\alpha\subset B$.  Then a $T$-chain beginning in $N_\alpha$ stays in $N_\alpha$ until leaving $B$ (at which point it also won't be in $N_{\alpha'}$ for any $\alpha'\neq \alpha$).  

The collection $\{N_\alpha\}$ covers $B$.  Now let $W_i^{\alpha} = U_i^{\alpha}\cap f^{-1}(N_{\alpha})$ and put $W_i = \bigcup_\alpha W_i^{\alpha}$.  Then $\{W_i\}_{i=0}^{k-1}$ is a cover of $f^{-1}(B)$ by $k$ open sets.  Since any $x\in f^{-1}(B)$ is in $f^{-1}(C_\alpha)$ for some $\alpha$, every $x\in f^{-1}(B)$ is contained in at least $k-m$ elements of $\{W_i\}_{i=0}^{k-1}$.  Moreover, since $f(W_i^{\alpha})\subset N_{\alpha}$, the properties of $N_\alpha$ imply that an $(F, T)$-component of $W_i$ is contained in an $F$-component of some $W_i^{\alpha}\subset U_i^{\alpha_j}$ and is hence $\mathcal{D}^{(k)}_f(F, S)$-bounded. \end{proof}


\lemma{(Finite union lemma, compare to \cite[Corollary 3.6]{Brodskiy2006AHT}) With the same setup as \ref{dad of a map}, let $A, B\subset X$, $F\in \mathcal{P}_{\text{fs}}(\Gamma)$ and $S\in \mathcal{P}_{\text{fs}}(\Lambda)$.  Assume that the $(F^{r}, S^{p})$-components of $A$ are $(F^{R}, S^{P})$-bounded, that the $(F^{l}, S^{q})$-components of $B$ are $(F^{L}, S^{Q})$-bounded; and that $2L + 2l<r$ and $2Q + 2q< p$.  Then the $(F^{l}, S^{q})$-components of $A\cup B$ are $(F^{r + R}, S^{p + P})\text{-bounded}$.}\label{finite union lemma}
\begin{proof}
Let $x_0, \ldots, x_n$ be an $(F^{l}, S^{q})$-chain in $A\cup B$.  Suppose $x_j$ and $x_k$ are two consecutive points contained only in $A$.  Then $x_{j+1}, \ldots, x_{k-1}$ is an $(F^{l}, S^{q})$-chain in $B$ and therefore is $(F^{L}, S^{Q})$-bounded.  

We therefore have that $x_{j+1}$ and $x_{k-1}$ are in the same $(F^{2L}, S^{2Q})$-component of $B$, so $x_i$ and $x_j$ are in the same $(F^{2L + 2l}, S^{2Q + 2q})$-component of $A$.  Since $2L + 2l < r$ and $2Q + 2q < p$, the points in the original chain which are contained only in $A$ therefore form an $(F^{r}, S^{p})$-chain in $A$ and so (considered as a subset of $A$) are $(F^{R}, S^{P})$-bounded.  The original chain, considered as a set, is therefore $(F^{2L + 2l + R}, S^{2Q + 2q + P})$-bounded.  This implies the $(F^{l}, S^{q})$-components of $A\cup B$ are $(F^{2L + 2l + R}, S^{2Q + 2q + P})$-bounded, which implies the result since $2L + 2l< r$ and $2Q + 2q<p$.  \end{proof}

\theorem{(compare to \cite[Theorem 4.9]{Brodskiy2006AHT}) Let $k = m + n + 1$, where $m, n\geq 0$.  Let $\Gamma\acts X$, $\Lambda\acts Y$ and $f: X\to Y$ be as in \ref{dad of a map}.  Suppose $\mathcal{D}_{Y}$ is an $n$-dimensional control function for $\Lambda\acts Y$ and $\mathcal{D}_f$ is an $m$-dimensional control function for $f$.  Then there is an $m+n$-dimensional control function for $\Gamma\acts X$ depending only on $\mathcal{D}_{Y}$, $\mathcal{D}_f$, $m$, and $n$.  } \label{control function theorem}

\begin{proof}Let $\alpha: \Gamma\times X\to \Lambda$ be the function described in \ref{dad of a map} associated to $f$.  Apply \ref{extra cover lemma} to $\mathcal{D}_{Y}$ to produce the function $\mathcal{D}^{(k)}_{Y}$.   Because we can always make the outputs of a control function larger sets (while keeping them finite) and it will still be a control function, we can assume that $\mathcal{D}_{Y}(T)\supset T$ for all $T\in \mathcal{P}_{\text{fs}}(\Lambda)$, so that $\mathcal{D}^{(k)}_{Y}(T)\supset T$ for all $T\in \mathcal{P}_{\text{fs}}(\Lambda)$.


Fix $F\in \mathcal{P}_{\text{fs}}(\Gamma)$.  Define inductively a sequence of finite subsets $S_Y^{(n+1)} \subset T_Y^{(n+1)}\subset S_Y^{(n)}\subset T_Y^{(n)}\subset \cdots \subset S_Y^{(1)}\subset T_Y^{(1)}\subset S_Y^{(0)}\subset T_Y^{(0)}$ starting from $S_Y^{(n+1)} = \alpha(F, X)$ and moving to lower indices so that for every $i$ we have $T_Y^{(i)} = \mathcal{D}^{(k)}_Y(S_Y^{(i)})$ and $S_Y^{(i)} = (T_Y^{(i+1)})^4$.  

Write $Y$ as the union of $n+1$ open sets $\{A_i\}_{i=1}^{n+1}$ such that all $S_Y^{(0)}$-components of $A_i$ are $T_Y^{(0)}$-bounded for every $i$.  Since $A_i\subset Y$, we can use \ref{extra cover lemma} to obtain a cover of $A_i$ by $k$ sets $\{U_i^j\}_{j=1}^k$ which are open in $X$ such that all $S_Y^{(i)}$-components of $U_i^j$ are $T_Y^{(i)}$-bounded for every $j$ and every point $y\in A_i$ belongs to at least $k-n = m+1$ sets.  To be clear, although $A_i$ is not necessarily closed, we can apply \ref{extra cover lemma} to $Y$ and then take intersections with $A_i$.  

Apply \ref{function extra cover lemma} to $\mathcal{D}_f$ and $F$ to obtain the function $\mathcal{D}^{(k)}_f$.     As with $\mathcal{D}_{Y}$, we can assume $\mathcal{D}_f(S, T)\supset S$ for all $S\in \mathcal{P}_{\text{fs}}(\Gamma)$ and $T\in \mathcal{P}_{\text{fs}}(\Lambda)$.  Define inductively another sequence of finite sets $S_X^{(n+1)} \subset T_X^{(n+1)}\subset S_X^{(n)}\subset T_X^{(n)}\subset \cdots \subset S_X^{(1)}\subset T_X^{(1)}$ starting with $S_X^{(n+1)} = F$ and for every $i$ defining $T_X^{(i)} = \mathcal{D}^{(k)}(S_X^{(i)}, T_Y^{(0)})$ and $S_X^{(i)} = (T_X^{(i+1)})^4$.  

By \ref{control function tricky lemma} and the properties of $A_i$, $f^{-1}(A_i)$ has a cover by $k$ sets $\{B_i^j\}_{j=1}^k$ which are open in $X$ such that all $(S_X^{(i)}, S_Y^{(0)})$-components are $(T_X^{(i)}, T_Y^{(0)})$-bounded for every $j$ and every point $x\in f^{-1}(A_i)$ belongs to at least $k-m = n + 1$ sets.  

Put $D_i^j = B_i^j\cap f^{-1}(U_i^j)$ and let $D^j$ be the union of all the $D_i^j$.  The $D^j$'s cover $X$ by Kolmogorov's argument: given $x\in X$, there is $i$ so that $f(x)\in A_i$.  The set of $j$'s such that $x\in B_i^j$ has at least $n + 1$ elements, and the set of $j$'s such that $f(x)\in U_i^j$ has at least $k-n$ elements, so these sets of such $j$'s cannot be disjoint by the pigeonhole principle.  

Notice all $(S_X^{(i)}, S_Y^{(i)})$-components of the set $D_i^j$ are $(T_X^{(i)}, T_Y^{(i)})$-bounded.  By \ref{finite union lemma} and induction (going from lower to higher indices), all $(S_X^{(n+1)}, S_Y^{(n+1)})$-components of the set $D^j$ are $((T_X^{(1)})^2, (T_Y^{(1)})^2)$-bounded.  Since $(S_X^{(n+1)}, S_Y^{(n+1)}) = (F, \alpha(F, X))$ by definition, and $(F, \alpha(F, X))$-components coincide with $F$-components, all $F$-components of $D^j$ are $(T_X^{(1)})^2$-bounded.  \end{proof}

\theorem{(Hurewicz-type theorem) Suppose $\Gamma\acts X$ and $\Lambda\acts Y$ are actions on compact metric spaces by homeomorphisms, and that $f: X\to Y$ and $\alpha: \Gamma\times X\to \Lambda$ are continuous functions such that $f(\gamma\cdot x) = \alpha(\gamma, x)\cdot f(x)$.  With $\dad(f)$ defined as in \ref{dad of a map}, it holds that $\dad_{free}(\Gamma\acts X)\leq \dad(f) + \dad_{free}(\Lambda\acts Y)$.  }\label{hurewicz-type theorem} 
\begin{proof}
Combine \ref{dimension implies function} with \ref{control function theorem}.  \end{proof}

\remark{The product theorem proven in the previous section can now be obtained as a consequence of \ref{hurewicz-type theorem}.  In fact, one can show a bit more by considering actions on bundles which respect the bundle structure.  The map from the bundle to its base space will have dynamic asymptotic dimension exactly that of the action on a fibre.  }

\corollary{Suppose $\Gamma\acts X$, $\Lambda\acts Y$, and $f: X\to Y$ are as above.  Then $f$ has an $m$-dimensional control function if and only if $\dad(f)\leq m$.  }\label{function iff dimension}
\begin{proof}
One direction was already proved in \ref{dimension implies function}.  For the other direction, apply \ref{control function theorem} to each subspace of $Y$ which is $0$-dimensional. \end{proof}

\section{Applications}\label{applications}


\normalfont

An important corollary of most Hurewicz-type theorems is that dimension is subadditive over extensions of objects.  Although there are natural statements of this kind that we can now prove for the dynamic asymptotic dimension, it is less obvious whether such theorems can be usefully applied.  Consider for instance, the action $\Z^2\acts S^1$ where each copy of $\Z$ acts by an irrational rotation.  The restriction to either copy of $\Z$ gives a minimal action so that the action of the other copy on the quotient by the orbit of the first is infinite dimensional (because the quotient is very non-Hausdorff).  Notwithstanding, we highlight two useful applications below.  The first shows that the dynamic asymptotic dimension does not change when passing to a finite index subgroup, which shows an additional aspect in which $\dad$ behaves like a large-scale dimension for group actions.  The second is an extension theorem for the dimension of odometer actions, which will have further applications in the final section.  

For some applications, it is easier to directly produce a control function rather than bound the dimension of a map.  In fact, this inspires an alternative definition for said dimension.  

\definition{Suppose $\Gamma\acts X$, $\Lambda\acts Y$, and $f: X\to Y$ are as in \ref{dad of a map}.  The value $\dad(f)$ can equivalently be defined to be the smallest $m$ such that $\dad(\{f^{-1}(\{y\}) : y\in Y\})\leq m$ uniformly in the sense of \ref{uniform dimension definition}, that is, $f^{-1}(\{y\})$ has an $m$-dimensional control function for each $y\in Y$, and these control functions can all be taken to be the same (independent of $y$).  } \label{new dad of a map}
\begin{proof}
If $\dad(f)\leq m$ according to the original definition, then $f$ has an $m$-dimensional control function by \ref{function iff dimension}.  This implies $\dad(\{f^{-1}(\{y\}) : y\in Y\})\leq m$ uniformly.  

Now assume $\dad(f^{-1}(\{y\}))\leq m$ uniformly.  Then using \ref{finite union lemma} and induction produces an $m$-dimensional control function for $f$.  \end{proof}

\corollary{Suppose $\Gamma\acts X$ is a free action on a compact metric space by homeomorphisms and that $\Lambda<\Gamma$ is a finite index subgroup.  Then $\dad(\Gamma\acts X) = \dad(\Lambda\acts X)$.  }
\begin{proof}
By passing to the normal core, we may assume without loss of generality that $\Lambda$ is normal.  Suppose $\dad(\Lambda\acts X) = m$.  Let $X/\Lambda$ be the quotient of $X$ by the orbit equivalence relation coming from $\Lambda\acts X$, and let $\Gamma/\Lambda\acts X/\Lambda$ be the obvious action.  Let $f: X\to X/\Lambda$ be the quotient map.  Since $\Gamma/\Lambda\acts X/\Lambda$ is locally finite, $\dad(\Gamma/\Lambda\acts X/\Lambda) = 0$.  Moreover, if $y\in X/\Lambda$, $f^{-1}(\{y\})$ is a $\Lambda$-orbit in $X$ and so, as $\Gamma\acts X$ is free, $\dad( f^{-1}(\{y\}))$ (with respect to $\Gamma\acts X$) is at most $\dad(\Lambda\acts f^{-1}(\{y\}))$, which in turn is at most $\dad(\Lambda\acts X)$.  Moreover, these inequalities hold quantitatively for all $y$ (i.e. if $\mathcal{D}$ is a control function for $\Lambda\acts X$, $F\mapsto \mathcal{D}(F\cap \Lambda)$ is a control function for $f^{-1}(\{y\})$ with respect to $\Gamma\acts X$).  Hence, $\dad(\Gamma\acts f^{-1}(\{y\}))\leq m$ uniformly, and so $\dad(f)\leq m$ by \ref{new dad of a map}.  The result then follows from \ref{hurewicz-type theorem}.  \end{proof}

\normalfont
%
%


For the extension theorem, we still need to be careful to avoid a similar pitfall as in \cite[Lemma 4.1(5)]{finnsell2015asymptotic} by assuming a certain compatibility between the extension and the filtration $(N_i)$.  Fortunately, we will see there are still many examples to which the theorem applies.  In fact, we will see there are examples where the extra assumption becomes immaterial.  

\definition{Suppose $\Gamma$ is a countable group, $\Delta<\Gamma$ is a subgroup, and $(N_i)$ is a countable collection of finite index normal subgroups of $\Gamma$ directed by inclusion.  We say $\Lambda$ is separated from $\Gamma$ by $(N_i)$ if for all $\gamma\in \Gamma\setminus \Delta$ there is $N_i$ such that $\gamma N_i\notin \{\delta N_i\mid \delta\in \Delta\}$.  Notice that in this case $\widehat{\Delta}_{(N_i\cap \Delta)}$ is closed as a subset of $\widehat{\Gamma}_{(N_i)}$.  }\label{separable subgroup definition}

\remark{If $\Gamma$ is a group, $1\to \Delta\to \Gamma\xrightarrow{q} \Lambda\to 1$ is an exact sequence, and $\Lambda$ is residually finite, there exists a family of finite index normal subgroups separating $\Delta$ from $\Gamma$.  For any sequence $(N_i)$ of finite-index normal subgroups of $\Gamma$, there is an associated sequence $(\tilde{N}_i)$ which separates $\Delta$ from $\Gamma$ and is a filtration if $N_i$ is.  }\label{separating example}
\begin{proof}
Let $\gamma\in \Gamma$ such that $\gamma\notin \Delta$.  Then $q(\gamma)$ is not the identity in $\Lambda$ and so there is $H_i$ such that $q(\gamma)H_i$ is not the identity.  This implies $(q^{-1}(H_i))$ separates $\Delta$ from $\Gamma$.  We can assume that $(H_i)$ is a sequence.  For the last claim, put $\tilde{N}_i = \cap_{j\leq i} (N_j\cap q^{-1}(H_j))$.  \end{proof}

%

\corollary{Suppose $\Gamma$ is a countable group and $(N_i)$ is a countable collection of finite index normal subgroups of $\Gamma$ directed by inclusion.  If $1\to \Delta\to \Gamma\to \Lambda\to 1$ is an extension of groups, $(N_i)$ induces collections of finite index subgroups $(N_i\cap \Delta)$ of $\Delta$ and $(N_i/(\Delta\cap N_i))$ of $\Lambda$.  If also $(N_i)$ separates $\Delta/(\cap_i N_i\cap \Delta)$ from $\Gamma/(\cap_i N_i)$, then $\dad_{free}(\Gamma\acts \widehat{\Gamma}_{(N_i)})\leq \dad_{free}(\Delta\acts \widehat{\Delta}_{(N_i\cap \Delta)}) + \dad_{free}(\Lambda\acts \widehat{\Lambda}_{(N_i/(\Delta\cap N_i))})$.  } \label{extension theorem}

\begin{proof}
Since we are using the definition of $\dad$ for free actions, any elements of $\Gamma$ which are trivial in every $N_i$ are immaterial, so we can assume $\Gamma\acts \widehat{\Gamma}_{(N_i)}$ is free and that $(N_i)$ separates $\Delta$ from $\Gamma$.  Supposing that $\dad(\Delta\acts \widehat{\Delta}_{(N_i\cap \Delta)})\leq m$ and $f: \widehat{\Gamma}_{N_i}\to \widehat{\Lambda}_{(N_i/(\Delta\cap N_i))}$ is obtained from the quotient maps $\Gamma/N_i\to \Lambda/(N_i/(\Delta\cap N_i))$, we will show $\dad(f)\leq m$ using \ref{new dad of a map}.  

Fix $F\in \mathcal{P}_{\text{fs}}(\Gamma)$.  Denote by $(e_i)$ the identity element in $\widehat{\Lambda}_{(N_i/\Delta\cap N_i)}$ and let $\iota: \widehat{\Delta}_{(N_i\cap \Delta)} \to \widehat{\Gamma}_{(N_i)}$ be the map coming from the inclusions $\Delta/(N_i\cap \Delta)\to \Gamma/N_i$.  Then $f^{-1}(\{(e_i)\}) = \iota(\widehat{\Delta}_{(N_i\cap \Delta)})$.  Since $\iota$ is equivariant and a homeomorphism onto its image, an $(m, F, S)$-cover of $\widehat{\Delta}_{(N_i\cap \Delta)}$ (which we have for some $S$ by assumption) gives rise to an $(m, F, S)$-cover of $\iota(\widehat{\Delta}_{(N_i\cap \Delta)})$ with respect to $\Delta\acts\widehat{\Gamma}_{(N_i)}$ (and vice versa).  Because $\Delta$ is separated from $\Gamma$ by $(N_i)$, any $\gamma\notin \Delta$ moves $x\in f^{-1}(\{(e_i)\})$ outside of $f^{-1}(\{(e_i)\})$.  Therefore, an $(m, F, S)$-cover of $f^{-1}(\{(e_i)\})$ with respect to $\Delta\acts\widehat{\Delta}_{(N_i\cap \Delta)}$ is also an $(m, F, S)$-cover of $f^{-1}(\{(e_i)\})$ with respect to $\Gamma\acts \widehat{\Gamma}_{(N_i)}$.  Since $f$ is a homomorphism of the groups $\widehat{\Gamma}_{(N_i)}$ and $\widehat{\Lambda}_{(N_i/\Delta\cap N_i)}$, the above argument can be repeated for any $x\in \widehat{\Lambda}_{(N_i/\Delta\cap N_i)}$ in place of $(e_i)$.  This follows since right multiplication by any element of $f^{-1}(\{x\})$ gives an orbit equivalence between $\Delta\acts f^{-1}(\{x\})$ and $\Delta\acts f^{-1}(\{(e_i)\})$.  Since $S$ does not depend on $x$, we are done.  \end{proof}


\corollary{Suppose $\Gamma$ is a countable group and $(N_i)$ is a sequence of finite index normal subgroups.  If $1\to \Delta\to \Gamma\to \Lambda\to 1$ is an extension of groups, $(N_i)$ induces collections of finite index subgroups $(N_i\cap \Delta)$ of $\Delta$ and $(N_i/(\Delta\cap N_i))$ of $\Lambda$.  If also $(N_i)$ separates $\Delta/(\cap_i N_i\cap \Delta)$ from $\Gamma/\cap_i N_i$, then $\asdim\Box_{(N_i)}\Gamma\leq \asdim\Box_{(N_i\cap \Delta)}\Delta + \asdim \Box_{(N_i/(\Delta\cap N_i))}\Lambda$.  } \label{extension theorem for box spaces}

\begin{proof} Combine \ref{extension theorem} with \ref{dimension of odometers and box spaces}.  \end{proof}

\section{Baumslag-Solitar groups} \label{baumslag-solitar groups section}
\normalfont
We conclude by demonstrating how the theoretical machinery introduced in this paper can be used to calculate the dimension of box spaces for many new examples.  What's more, we would like to prove special cases of the (unproved) theorems in \cite[Section 4]{finnsell2015asymptotic} about elementary amenable groups.  The idea is to use \cite{MR1143191}*{page 238}, which describes in particular how any elementary amenable group can be constructed by iterated extensions $1\to \Delta\to \Gamma\to \Lambda\to 1$ where $\Lambda$ is virtually abelian; together with an extension theorem for the asymptotic dimension of box spaces.  Unfortunately, the requirements that $(N_i)$ be a sequence with certain separation properties creates problems when trying to replicate the proofs of these results.  It should be noted however that we are `close' to being able to make these arguments work, and that our consideration of separable subgroups is a new insight which may help in implementing the approach attempted in \cite{finnsell2015asymptotic}.  Further to the point, we observe the following: since the group $\Lambda$ is virtually abelian and hence residually finite, there exists for any such extension a family of finite index normal subgroups of $\Gamma$ separating $\Delta$ from $\Gamma$.  As \ref{direct limit lemma} lets us handle direct limit groups, it seems plausible that one could prove an extension theorem for box spaces directly (keeping the assumption that $(N_i)$ separates $\Delta$ from $\Gamma$, and such an extension theorem would not require $(N_i)$ to be a sequence.  In this way, one could obtain results similar to the statements in \cite[Section 4]{finnsell2015asymptotic} using full box spaces, but we leave that for another time and instead focus on a few interesting examples.  

We will show that for all integers $n$, the box spaces (coming from sequential filtrations) of the Baumslag-Solitar groups $\text{BS}(1, n)$ have dimension $2$.  This class includes infinitely-many groups which are not virtually nilpotent (in fact, which have exponential growth), providing concrete evidence for the conjecture that finite-dimensional amenable groups should have finite-dimensional box spaces.    

\definition{The Baumslag-Solitar group $\text{BS}(m, n)$ is the group given by the presentation $\langle a, b \mid ba^mb^{-1} = a^n\rangle$.  }

\lemma{Suppose $\Gamma = \text{BS}(1, n)$ for some $n$ and $N\subset \Gamma$ is a normal subgroup of finite index.  Then $N = \langle a^l, b^k\rangle$ for some $l, k\geq 1$ and $l\mid nk-1$ (that is $l$ divides $nk-1$).  }

\begin{proof}
If $N\subset \langle a\rangle$ then $\Gamma/N$ is not finite, so $N$ contains some power of $b$.  Let $k\geq 1$ be the least integer such that $b^k\in N$.  Since $N$ is normal in $\Gamma$, $a^{-1}b^ka = a^{nk-1}b^k\in N$ and so $a^{nk-1}\in N$.  There is therefore some $l\geq 1$ dividing $nk-1$ such that $a^l\in N$.  Let $l\geq 1$ be the least integer with these properties.  Then $N = \langle a^l, b^k\rangle$.  \end{proof}

\lemma{Suppose $\Gamma = \text{BS}(1, n)$.  Then $\Gamma$ fits into the exact sequence $1\to \Z\to \Gamma\to \Z\to 1$.  Moreover, if $(N_i)$ is a filtration of $\Gamma$, then $(N_i)$ separates $\Z$ from $\Gamma$.  } \label{separating example lemma}

\begin{proof}
The first map in the sequence can be thought of as the inclusion of $\langle a\rangle$ into $\Gamma$.  The second is then the quotient by the normal subgroup $\langle a\rangle$.  

Suppose $(N_i)$ is a filtration of $\Gamma$ by normal subgroups.  Let $\gamma\in \Gamma$ with $\gamma\notin \langle a\rangle$.  Any element of $\Gamma$ can be written as $a^mb^j$ for some $m, j\in \Z$, so $\gamma = a^mb^j$ with $j\neq 0$.  By the previous lemma, each $N_i$ has the form $\langle a^{l_i}, b^{k_i}\rangle$ for $l_i, k_i\geq 1$.  Since $N_i$ is a filtration, we have $i$ such that $k_i$ does not divide $j$, since otherwise we would have $b^j\in N_i$ for all $i$.  Let $\pi: \Gamma\to \Gamma/N_i$ be the quotient map.  The quotient group $\Gamma/N_i$ has the presentation $\langle a, b\mid a^{l_i} = b^{k_i} = 1, bab^{-1} = a^n\rangle$.  Using the same normal form for elements, we see that $b^{j} = \pi(b)^j\notin \langle a\rangle = \langle \pi(a)\rangle$.  \end{proof}

\normalfont
We will calculate the asymptotic dimension of $\text{BS}(1, n)$ for all $n$ and use that to give the lower bound for the dimension of its box space.  To do this, we will combine some results already in the literature to give an affirmative answer to \cite[Question 1.4]{asdimnote}.  More precisely, \cite[Theorem 2]{MR2431020} says that a finitely-presented group of asymptotic dimension $1$ is virtually free, and \cite[Theorem 3]{MR223439} says that a finitely generated, virtually free, torsion free group is free.  Combining these results therefore yields the following.  


\theorem{A finitely-presented, torsion free group of asymptotic dimension $1$ is free.  }\label{asdim free}
%
%

\normalfont Of course, by \ref{dimension of odometers and box spaces}, the box spaces in the following examples can be replaced by the corresponding odometers, and the statements still hold with $\dad$ in place of asymptotic dimension.  

\corollary{Let $\Gamma = \text{BS}(1, n)$.  Then $\Gamma$ has exponential growth (so in particular is not virtually nilpotent) for $n\geq3$ and has $\asdim\Box_{(N_i)}\Gamma = \asdim \Gamma = 2$ for any sequence $(N_i)$ of normal subgroups which is a filtration of $\Gamma$.  } \label{example}

\begin{proof}
The group $\Gamma$ is residually finite by \cite[Theorem C]{MR285589} and fits into an exact sequence $1\to \Z\to \Gamma \to \Z\to 1$, so \cite[Theorem 68]{asdimBD} (subadditivity over extensions) shows that $\asdim\Gamma \leq 2$.  Because $\Gamma$ is torsion-free but not free, \ref{asdim free} implies it cannot be that $\asdim\Gamma = 1$, and so $\asdim\Gamma = 2$.  

By \ref{separating example lemma}, $(N_i)$ also separates $\Z$ from $\Gamma$ in this extension, and so \ref{extension theorem for box spaces} implies $\asdim\Box_{(N_i)}\Gamma\leq 2$.  The reverse inequality follows from \cite[Proposition 3.1]{boxspacesDT} since $\asdim\Gamma = 2$ and $(N_i)$ is a filtration.  \end{proof}




\normalfont
We can also generalize the above examples by defining more complicated groups with a similar presentation.  An upper bound for the dimension of box spaces of such groups can be found in much the same way.  While we need to assume the collection $(N_i)$ has special separating properties, \ref{separating example} shows there are examples of such collections.  

\corollary{Let $m_1, m_2\geq 1$ be integers, $n_{i, j}\in \Z$ for all $1\leq i\leq m_1$ and $1\leq j\leq m_2$, and $\sigma_j\in S_{m_1}$ for $1\leq j\leq m_2$ (here $S_{m_1}$ denotes permutations of $\{1, \ldots, m_1\}$).  Now define 

$$\Gamma := \langle a_1, \ldots, a_{m_1}, b_1, \ldots, b_{m_2} \mid [a_i, a_j] = 1, \text{ } [b_i, b_j] = 1, \text{ } b_ja_ib_j^{-1} = a_{\sigma_j(i)}^{n_{i, j}}\rangle.$$

\noindent Then $\asdim \Box_{(N_i)} \Gamma\leq m_1 + m_2$ for any sequence $(N_i)$ of finite index normal subgroups which separates $\langle a_1, \ldots, a_{m_1}\rangle$ from $\Gamma$ (and such a sequence exists).  }

\begin{proof}

The group $\Gamma$ fits into the exact sequence $1\to \Z^{m_1}\to \Gamma\to \Z^{m_2}\to 1$.  To be precise, the first map identifies $\Z^{m_1}$ with the subgroup generated by the $a_i$'s, and this subgroup is moreover normal in $\Gamma$.  Therefore, $\asdim\Box_{(N_i)} \Gamma\leq m_1 + m_2$.  As mentioned earlier, a sequence which separates $\Z^{m_1}$ from $\Gamma$ exists by \ref{separating example}.  \end{proof}
 

\normalfont There are also groups known as generalized Baumslag-Solitar groups \cite{DELGADO2011398}.  Most of these are non-amenable, and so will have infinite dimensional box spaces, but their consideration leads to at least one more interesting example.  

\corollary{Let $\Gamma = \langle a, b \mid a^2 = b^{2}\rangle$.  Then $\asdim \Box_{(N_i)}\Gamma = 2$ for any sequence $(N_i)$ of finite index normal subgroups which separates $\langle a^2\rangle$ from $\Gamma$ (and such a sequence exists).  }

\begin{proof}
Once again, \ref{separating example} shows we have a sequence $(N_i)$ satisfying the hypotheses.  For the upper bound, observe that $\Gamma$ fits into the exact sequence $1\to \Z\to \Gamma\to \Z/2*\Z/2\to 1$ where the first map is the inclusion of the normal subgroup $\langle a^2\rangle$ and the second map is the quotient by this subgroup.  Since $\Z/2*\Z/2$ is the infinite dihedral group (which is virtually cyclic), $\asdim \Box_{(N_i)}\Gamma \leq 2$.  
 
Since $\Gamma/\langle a^{2n}\rangle\cong \Z/2n*\Z/2n$ is residually finite for all $n$, and any non-trivial element $\gamma\in \Gamma$ is non-trivial in $\Gamma/\langle a^{2n}\rangle$ for some $n$, it follows that $\Gamma$ is residually finite.  Hence, $\asdim\Box_{(N_i)}\Gamma\geq \asdim\Gamma$ by \cite[Proposition 3.1]{boxspacesDT}, and we are done if we show $\asdim \Gamma \geq 2$.  Since $\Gamma$ is torsion free but not free, this follows from \ref{asdim free}. \end{proof}

\begin{acknowledgements}
\noindent \raggedright This work was partly supported by the NSF (DMS 1901522) and has received funding from the European Research Council (ERC) under the European Union's Horizon 2020 research and innovation programme.  
\end{acknowledgements}

\bibliography{mybibliography2.bib}
\bibliographystyle{plain}

\end{document}